\documentclass{amsart}
\usepackage{amsmath}
\usepackage{amsfonts,amssymb,amsthm}
\usepackage{mathscinet}
\makeatletter
\newtheorem{theorem}{Theorem}[section]
\newtheorem{lemma}[theorem]{Lemma}
\newtheorem{proposition}[theorem]{Proposition}    
\newtheorem{corollary}[theorem]{Corollary}
\theoremstyle{definition}
\newtheorem{definition}[theorem]{Definition}
\newtheorem{remark}[theorem]{Remark}
\newtheorem{example}[theorem]{Example}

\begin{document}
\title[B-J orthogonality preservers between non-isometric spaces]{On Birkhoff-James orthogonality preservers between real non-isometric Banach spaces}
\author{Ryotaro Tanaka}
\address{Katsushika Division, Institute of Arts and Sciences, Tokyo University of Science, Tokyo 125-8585, Japan}
\email{r-tanaka@rs.tus.ac.jp}
\date{}
\begin{abstract}
Real smooth three-dimensional or higher Banach spaces are isomorphic with respect to the nonlinear structure of Birkhoff-James orthogonality if and only if they are isometrically isomorphic. Moreover, using smooth Radon planes and non-smooth direct sums, in arbitrary dimensions, we construct examples of non-isometric pairs of real Banach spaces that admit norm-preserving homogeneous bicontinuous Birkhoff-James orthogonality preservers among them.
\end{abstract}

\maketitle
%%%%%%%%%%%%%%%%%%%%%%%%%%%%%%%%%%%%%%%%%%%%%%%%%%%%%%%%%%
%%%%%%%%%%%%%%%%%%%%%%%%%%%%%%%%%%%%%%%%%%%%%%%%%%%%%%%%%%
\section{Introduction}

Throughout this paper, the term ``Banach space'' indicates a \emph{real} Banach space. The present paper is concerned with nonlinear classification of Banach spaces based on Birkhoff-James orthogonality. Let $X$ be a Banach space, and let $x,y \in X$. Then, $x$ is said to be \emph{Birkhoff-James} orthogonal to $y$, which is denoted as $x \perp_{BJ}y$, if $\|x+\lambda y\| \geq \|x\|$ for each $\lambda \in \mathbb{R}$. This generalized orthogonality relation in Banach spaces was introduced by Birkhoff~\cite{Bir35} and deeply studied by James~\cite{Jam47a,Jam47b}, and is known as one of the most important generalized orthogonality relations since it is closely related to the geometric structure of Banach spaces. For example, if $x,y \in X$ and $x \perp_{BJ}y$, then $x$ is \emph{a} nearest point from the origin in the line $L=\{ x+\lambda y :\lambda \in \mathbb{K}\}$. Further, if $x$ is a unit vector, $L$ becomes a tangent line to the unit ball of $X$. This last statement also implies a relationship between Birkhoff-James orthogonality and support functionals for the unit ball $B_X$ of $X$. Indeed, studies on Birkhoff-James orthogonality have frequently begun with recalling the following useful result by James~\cite{Jam47b}:
%%%%%%%%%%%%%%%%%%%%%%%%%%%%%%%%%%%%%%%%%%%%%%%%%%%%%%%%%%
\begin{lemma}[James, 1947]\label{James}

Let $X$ be a Banach space, and let $x,y \in X$. Then, $x \perp_{BJ}y$ if and only if $f(y)=0$ for some $f \in \nu (x)$, where $\nu (x)$ is the set of all elements of $B_{X^*}$ that support $B_X$ at $x$; that is, $\nu (x)=\{ f\in B_{X^*} : f(x)=\|x\|\}$.
\end{lemma}
%%%%%%%%%%%%%%%%%%%%%%%%%%%%%%%%%%%%%%%%%%%%%%%%%%%%%%%%%%
It is known that Birkhoff-James orthogonality is \emph{non-degenerate} (that is, $x \perp_{BJ} x$ implies $x=0$), and \emph{homogeneous} (that is, $x \perp_{BJ} y$ implies $\alpha x\perp_{BJ} \beta y$ for each $\alpha ,\beta \in \mathbb{R}$). However, it is \emph{asymmetric} in general. Indeed, as was independently shown by Day~\cite[Theorem 6.4]{Day47} and James~\cite[Theorem 1]{Jam47a}, if $\dim X \geq 3$ and $x \perp_{BJ} y$ implies $y \perp_{BJ} x$ in $X$, then $X$ is a Hilbert space. Here, it should be noted that the assumption on the dimension of $X$ cannot be omitted. We can construct non-Hilbert two-dimensional real Banach spaces in which Birkhoff-James orthogonality is symmetric (called \emph{Radon planes}); see, for example,~\cite{Bir35} and \cite{Day47}, and~\cite{KST20} and \cite{MS06}, for various characterizations of Radon planes. For further information about generalized orthogonality relations, the readers are referred to Amir~\cite{Ami86}, who described many characterizations of inner product spaces, and a comprehensive survey on Birkhoff-James orthogonality and isosceles orthogonality by Alonso, Martini, and Wu~\cite{AMW12}.

As an interesting application of Birkhoff-James orthogonality to the theory of classification of Banach spaces, Koldobsky~\cite{Kol93} showed the following theorem.
%%%%%%%%%%%%%%%%%%%%%%%%%%%%%%%%%%%%%%%%%%%%%%%%%%%%%%%%%%
\begin{theorem}[Koldobsky, 1993]\label{BT-main}

Let $X,Y$ be Banach spaces, and let $T:X \to Y$ be linear. Then, $T$ is a scalar multiple of a linear isometry if and only if $x,y \in X$ and $x \perp_{BJ}y$ imply $Tx \perp_{BJ}Ty$.
\end{theorem}
%%%%%%%%%%%%%%%%%%%%%%%%%%%%%%%%%%%%%%%%%%%%%%%%%%%%%%%%%%
Later, the complex version of the preceding theorem was proved by Blanco and Turn\v{s}ek~\cite{BT06}. Moreover, in 2019, W\'{o}jcik~\cite{Woj19} gave a simpler proof of Theorem~\ref{BT-main}, and improved it as follows.
%%%%%%%%%%%%%%%%%%%%%%%%%%%%%%%%%%%%%%%%%%%%%%%%%%%%%%%%%%
\begin{theorem}[W\'{o}jcik, 2019]\label{W-main}

Let $X,Y$ be Banach spaces, and let $T:X \to Y$ be additive. Then, $T$ is a scalar multiple of a linear isometry if and only if $x,y \in X$ and $x \perp_{BJ}y$ imply $Tx \perp_{BJ}Ty$.
\end{theorem}
%%%%%%%%%%%%%%%%%%%%%%%%%%%%%%%%%%%%%%%%%%%%%%%%%%%%%%%%%%
An important consequence of Theorem~\ref{W-main} is the following: For Banach spaces $X,Y$, if there exists an additive bijection $T:X \to Y$ that preserves Birkhoff-James orthogonality in \emph{one direction}, then $X=Y$ (that is, $X$ is \emph{isometrically isomorphic} to $Y$). This means that the combination of structures of addition and Birkhoff-James orthogonality determines the entire structure of a Banach space. In other words, with the aid of additivity, Banach spaces are fully classified by their structure of Birkhoff-James orthogonality. A natural question then arises: What happens if the additivity is omitted? In this direction, the author~\cite{Tan22} studied the nonlinear equivalence of Banach spaces based on Birkhoff-James orthogonality. Let $X,Y$ be Banach spaces, and let $T:X \to Y$ be a (possibly non-additive) bijection. Then, $T$ is called a \emph{Birkhoff-James orthogonality preserver} if it preserves the Birkhoff-James orthogonality in \emph{both directions} (that is, $x \perp_{BJ} y$ if and only if $Tx \perp_{BJ} Ty$); in addition, $X$ is said to be isomorphic to $Y$ \emph{with respect to the structure of Birkhoff-James orthogonality}, denoted by $X\sim_{BJ}Y$, a Birkhoff-James orthogonality preserver exists between $X$ and $Y$. Here, the bidirectionality of a Birkhoff-James orthogonality preserver is required to ensure the symmetry of ``$\sim_{BJ}$,'' whereas the bijectivity is assumed to  guarantee that $T(X)$ is also a Banach space. It was shown in \cite[Corollary 3.7 and Theorem 3.16]{Tan22} that if $X,Y$ are finite-dimensional, or reflexive and smooth, then $X\sim_{BJ}Y$ implies $X \cong Y$ (that is, $X$ is \emph{linearly isomorphic} to $Y$). Moreover, \cite[Theorem 4.3]{Tan22} states that if $X$ is a Banach space with $\dim X \geq 3$, and $H$ is a Hilbert space, then $X\sim_{BJ}H$ implies $X=H$. By contrast, an example of a non-Hilbert smooth Radon plane $X$ such that $X \sim_{BJ}\ell_2^2$ was constructed in \cite[Example 4.8]{Tan22}, where $\ell_2^2$ denotes the two-dimensional Euclidean space. This last result shows that, at least in the two-dimensional setting, a Birkhoff-James orthogonality preserver is not necessarily a scalar multiple of an isometric isomorphism. Hence, in general, the additivity assumption cannot be omitted in Theorem~\ref{W-main}.

In this paper, we present two theorems that further develop the theory of nonlinear classification of Banach spaces through Birkhoff-James orthogonality. The first one states that if $X,Y$ are smooth Banach spaces, then $X \sim_{BJ}Y$ if and only if $X=Y$. The second gives non-isometric pairs of Banach spaces of arbitrary dimensions that admit norm-preserving homogeneous bicontinuous Birkhoff-James orthogonality preservers between them.
%%%%%%%%%%%%%%%%%%%%%%%%%%%%%%%%%%%%%%%%%%%%%%%%%%%%%%%%%%

%%%%%%%%%%%%%%%%%%%%%%%%%%%%%%%%%%%%%%%%%%%%%%%%%%%%%%%%%%
%%%%%%%%%%%%%%%%%%%%%%%%%%%%%%%%%%%%%%%%%%%%%%%%%%%%%%%%%%

%%%%%%%%%%%%%%%%%%%%%%%%%%%%%%%%%%%%%%%%%%%%%%%%%%%%%%%%%%
%%%%%%%%%%%%%%%%%%%%%%%%%%%%%%%%%%%%%%%%%%%%%%%%%%%%%%%%%%
\section{Results in smooth Banach spaces}

Here in, we show that Birkhoff-James orthogonality preservers between non-isometric Banach spaces cannot be constructed under a smooth setting, where a Banach space $X$ is said to be \emph{smooth} if $\nu (x)$ is a singleton for each nonzero $x \in X$. To see this, we recall the following results described in~\cite[Theorems 3.6 and 3.9]{Tan22}.
%%%%%%%%%%%%%%%%%%%%%%%%%%%%%%%%%%%%%%%%%%%%%%%%%%%%%%%%%%
\begin{theorem}\label{dim-eq}

Let $X,Y$ be Banach spaces. Suppose that either $X$ or $Y$ is finite-dimensional. Then, $X \sim_{BJ}Y$ implies that $\dim X=\dim Y$
\end{theorem}
%%%%%%%%%%%%%%%%%%%%%%%%%%%%%%%%%%%%%%%%%%%%%%%%%%%%%%%%%%
\begin{theorem}\label{sm-ref}

Let $X,Y$ be smooth Banach spaces, and let $T:X \to Y$ be a Birkhoff-James orthogonality preserver. If $M$ is a reflexive subspace of $X$, then $T(M)$ is a closed subspace of $Y$.
\end{theorem}
%%%%%%%%%%%%%%%%%%%%%%%%%%%%%%%%%%%%%%%%%%%%%%%%%%%%%%%%%%
The argument herein is essentially based on the fundamental theorem of projective geometry. To be precise, we make use of the following form given by Mackey~\cite[Lemma A]{Mac42}.
%%%%%%%%%%%%%%%%%%%%%%%%%%%%%%%%%%%%%%%%%%%%%%%%%%%%%%%%%%
\begin{lemma}[Mackey, 1942]\label{Mackey}

Let $X,Y$ be Banach spaces, and let $\mathcal{M}_1(X)$ and $\mathcal{M}_1(Y)$ be the families of all one-dimensional subspaces of $X$ and $Y$, respectively. If $\rho :\mathcal{M}_1(X) \to \mathcal{M}_1(Y)$ is a bijection that preserves the linear independence of the finite elements of $\mathcal{M}_1(X)$ and $\mathcal{M}_1(Y)$. Then, there exists a linear bijection $T:X \to Y$ such that $T(M)=\rho (M)$ for each $M \in \mathcal{M}_1(X)$.
\end{lemma}
%%%%%%%%%%%%%%%%%%%%%%%%%%%%%%%%%%%%%%%%%%%%%%%%%%%%%%%%%%
Furthermore, the following auxiliary result is needed. In the following, if $x_1,\ldots ,x_n$ are elements of a Banach space, then $[x_1,\ldots ,x_n]$ denotes the linear span of $x_1,\ldots ,x_n$.

%%%%%%%%%%%%%%%%%%%%%%%%%%%%%%%%%%%%%%%%%%%%%%%%%%%%%%%%%%
\begin{lemma}\label{order-iso-one}

Let $X,Y$ be Banach spaces, and let $\mathcal{F}(X)$ and $\mathcal{F}(Y)$ be the families of all finite-dimensional subspaces of $X$ and $Y$, respectively. If $\rho:\mathcal{F}(X) \to \mathcal{F}(Y)$ is an order isomorphism, then $\rho$ is restricted to a bijection from $\mathcal{M}_1(X)$ onto $\mathcal{M}_1(Y)$ that preserves the linear independence of the finite elements of $\mathcal{M}_1(X)$ and $\mathcal{M}_1(Y)$.
\begin{proof}
We first note that $\rho (\{0_X\})=\{0_Y\}$. Indeed, we have $\{0_Y\}=\rho (M)$ for some $M \in \mathcal{F}(X)$. Since $\{0_X\} \subset M$, it follows that $\rho (\{0_X\}) \subset \rho (M)=\{0_Y\}$. This proves that $\rho (\{0_X\})=\{0_Y\}$.

Next, suppose that $x \in X \setminus \{0\}$. Then, $\rho ([x])\neq \{0_Y\}$ as indecated in the preceding paragraph. Take an arbitrary $y \in \rho ([x])\setminus \{0\}$. It follows from $[y] \subset \rho ([x])$ that $\rho^{-1}([y]) \subset [x]$; by $\rho^{-1}([y])\neq \{0_X\}$, it implies that $\rho^{-1}([y])=[x]$. Hence, we obtain $\rho ([x])=[y]$. In particular, $\rho (\mathcal{M}_1(X)) \subset \mathcal{M}_1(Y)$. Since $\rho^{-1}$ has the same property as $\rho$, we also have $\rho^{-1}(\mathcal{M}_1(Y)) \subset \mathcal{M}_1(X)$. Therefore, $\rho$ is restricted to a bijection from $\mathcal{M}_1(X)$ onto $\mathcal{M}_1(Y)$.

Now, let $x_1,\ldots ,x_n \in X$ occur such that $\{x_1,\ldots ,x_n\}$ is linearly independent. For each $j$, choose a $y_j \in \rho ([x_j]) \setminus \{0\}$. Then, $\rho ([x_j])=[y_j]$ as indicated in the preceding paragraph. If $\{y_1,\ldots ,y_n\}$ is linearly dependent, we have $[y_i] \subset [y_1,\ldots ,y_{i-1},y_{i+1},\ldots ,y_n]$ for some $i$. Meanwhile, if $j \neq i$, it follows from $[x_j] \subset [x_1,\ldots ,x_{i-1},x_{i+1},\ldots ,x_n]$ that $[y_j] \subset \rho ([x_1,\ldots ,x_{i-1},x_{i+1},\ldots ,x_n])$, which implies that
\[
[y_i] \subset [y_1,\ldots ,y_{i-1},y_{i+1},\ldots ,y_n] \subset \rho ([x_1,\ldots ,x_{i-1},x_{i+1},\ldots ,x_n]) .
\]
Thus, we obtain $[x_i] \subset [x_1,\ldots ,x_{i-1},x_{i+1},\ldots ,x_n]$. This contradicts the linear independence of $\{x_1,\ldots ,x_n\}$. Hence, $\{y_1,\ldots ,y_n\}$ must be linearly dependent. The exact same argument is also valid for $\rho^{-1}$. Consequently, $\rho$ preserves the linear independence of the finite elements of $\mathcal{M}_1(X)$ and $\mathcal{M}_1(Y)$.
\end{proof}
\end{lemma}
%%%%%%%%%%%%%%%%%%%%%%%%%%%%%%%%%%%%%%%%%%%%%%%%%%%%%%%%%%
Combining these results with Theorem~\ref{BT-main}, we have an improvement of \cite[Theorem 3.16]{Tan22}. The following is the first main theorem of this paper.
%%%%%%%%%%%%%%%%%%%%%%%%%%%%%%%%%%%%%%%%%%%%%%%%%%%%%%%%%%
\begin{theorem}\label{sm-eq}

Let $X,Y$ be smooth Banach spaces. Suppose that $\dim X \geq 3$ or $\dim Y \geq 3$. Then, $X\sim_{BJ}Y$ if and only if $X=Y$.
\begin{proof}
It is sufficient to show the ``only if'' part. Suppose that $X\sim_{BJ}Y$. Then, there exists a Birkhoff-James orthogonality preserver $T:X \to Y$. We first note that $\dim X \geq 3$ and $\dim Y \geq 3$ by Theorem~\ref{dim-eq}. Let $M$ be a finite-dimensional subspace of $X$. Since $M$ is reflexive, by Theorem~\ref{sm-ref}, $T(M)$ is a closed subspace of $Y$. In particular, the restriction of $T$ to $M$ gives a Birkhoff-James orthogonality preserver from $M$ onto $T(M)$. Hence, $\dim M=\dim T(M)$ by Theorem~\ref{dim-eq}. This proves that $T$ maps a finite-dimensional subspace of $X$ to a subspace of $Y$ having the same dimension. Since $T^{-1}$ has the same property as $T$, we can define an order isomorphism $\rho :\mathcal{F}(X) \to \mathcal{F}(Y)$ by $\rho (M)=T(M)$. Moreover, by Lemma~\ref{order-iso-one}, $\rho$ is restricted to a bijection from $\mathcal{M}_1(X)$  onto $\mathcal{M}_1(Y)$ that preserves the linear independence of the finite elements of $\mathcal{M}_1(X)$ and $\mathcal{M}_1(Y)$. From this, by Lemma~\ref{Mackey}, we have a linear bijection $S:X \to Y$ such that $S(M)=T(M)$ for each $M \in \mathcal{M}_1(X)$. Now, suppose that $x,y \in X$, and $x\perp_{BJ}y$. Then, $Tx \perp_{BJ} Ty$. Based on the homogeneity of Birkhoff-James orthogonality, and $[Sx]=[Tx]$ and $[Sy]=[Ty]$, it is implied that $Sx \perp_{BJ} Sy$. Thus, by Theorem~\ref{BT-main}, $S$ is a scalar multiple of an isometric isomorphism from $X$ onto $Y$. This proves that $X=Y$.
\end{proof}
\end{theorem}
%%%%%%%%%%%%%%%%%%%%%%%%%%%%%%%%%%%%%%%%%%%%%%%%%%%%%%%%%%
\begin{remark}
In contrast to the previous version of the preceding theorem~\cite[Theorem 3.16]{Tan22}, the proof of Theorem~\ref{sm-eq} no longer relies on the main results of Mackey~\cite[Theorem in p. 246]{Mac42}. The property of the existing linear operator $S$ follows from Theorem~\ref{BT-main}.
\end{remark}
%%%%%%%%%%%%%%%%%%%%%%%%%%%%%%%%%%%%%%%%%%%%%%%%%%%%%%%%%%
\begin{remark}

The assumption that $\dim X \geq 3$ or $\dim Y \geq 3$ cannot be omitted. Indeed, it was noted in~\cite[Example 4.8]{Tan22} that a two-dimensional non-Hilbert smooth real Banach space $X$ exists that satisfies $X\sim_{BJ}\ell_2^2$.
\end{remark}
%%%%%%%%%%%%%%%%%%%%%%%%%%%%%%%%%%%%%%%%%%%%%%%%%%%%%%%%%%

%%%%%%%%%%%%%%%%%%%%%%%%%%%%%%%%%%%%%%%%%%%%%%%%%%%%%%%%%%
%%%%%%%%%%%%%%%%%%%%%%%%%%%%%%%%%%%%%%%%%%%%%%%%%%%%%%%%%%

%%%%%%%%%%%%%%%%%%%%%%%%%%%%%%%%%%%%%%%%%%%%%%%%%%%%%%%%%%
%%%%%%%%%%%%%%%%%%%%%%%%%%%%%%%%%%%%%%%%%%%%%%%%%%%%%%%%%%
\section{Tools for construction}

To prove the second main theorem of this paper, preliminary works are needed. We first recall the acute and obtuse angles based on Birkhoff-James orthogonality.
%%%%%%%%%%%%%%%%%%%%%%%%%%%%%%%%%%%%%%%%%%%%%%%%%%%%%%%%%%
\begin{definition}

Let $X$ be a Banach space, and let $x,y \in X$. Then, $x$ is \emph{at an acute angle} to $y$, as denoted by $x \perp_{BJ}^+ y$, if $\|x+\lambda y\| \geq \|x\|$ for each $\lambda \geq 0$. If $x \perp_{BJ}^+ y$ and $x \not \perp_{BJ} y$, then $x$ is \emph{at a strictly acute angle} to $y$, which we write as $x \perp_{BJ}^{++}y$.
\end{definition}
%%%%%%%%%%%%%%%%%%%%%%%%%%%%%%%%%%%%%%%%%%%%%%%%%%%%%%%%%%
\begin{definition}

Let $X$ be a Banach space, and let $x,y \in X$. Then, $x$ is \emph{at an obtuse angle} to $y$, as denoted by $x \perp_{BJ}^- y$, if $\|x+\lambda y\| \geq \|x\|$ for each $\lambda \leq 0$. If $x \perp_{BJ}^- y$ and $x \not \perp_{BJ} y$, then $x$ is \emph{at a strictly obtuse angle} to $y$, which we write as $x \perp_{BJ}^{--}y$.
\end{definition}
%%%%%%%%%%%%%%%%%%%%%%%%%%%%%%%%%%%%%%%%%%%%%%%%%%%%%%%%%%
\begin{remark}

The notions of Birkhoff-James acute and obtuse angles were already considered by Sain~\cite{Sai17}. In particular, their basic properties were stated in \cite[Proposition 2.1]{Sai17}. For example, we have the following properties of ``$\perp_{BJ}^+$'' and ``$\perp_{BJ}^-$'':
\begin{itemize}
\item[(i)] $x \perp_{BJ} y$ if and only if $x \perp_{BJ}^+y$ and $x \perp_{BJ}^-y$,
\item[(ii)] $x \perp_{BJ}^+ y$ implies $\alpha x \perp_{BJ}^+ \beta y$ for each nonnegative real numbers $\alpha,\beta$,
\item[(iii)] $x \perp_{BJ}^- y$ implies $\alpha x \perp_{BJ}^- \beta y$ for each nonnegative real numbers $\alpha,\beta$, and
\item[(iv)] $x\perp_{BJ}^+y$ if and only if $x \perp_{BJ}^- (-y)$.
\end{itemize}
\end{remark}
%%%%%%%%%%%%%%%%%%%%%%%%%%%%%%%%%%%%%%%%%%%%%%%%%%%%%%%%%%
Now, using support functionals, we present the following characterizations of Birkhoff-James acute and strictly acute angles. An idea is derived from Lemma~\ref{James}.

%%%%%%%%%%%%%%%%%%%%%%%%%%%%%%%%%%%%%%%%%%%%%%%%%%%%%%%%%%
\begin{lemma}\label{acute}

Let $X$ be a Banach space, and let $x,y \in X$. Then, the following hold:
\begin{itemize}
\item[{\rm (i)}] $x \perp_{BJ}^+y$ if and only if $f(y) \geq 0$ for some $f \in \nu (x)$,
\item[{\rm (ii)}] $x \perp_{BJ}^{++}y$ if and only if $f(y)>0$ for each $f \in \nu (x)$.
\end{itemize}
\end{lemma}
\begin{proof}
(i) Suppose that $x \perp_{BJ}^+y$. Then, for each $n \in \mathbb{N}$, there exists an $f_n \in S_{X^*}$ such that
\[
f_n(x+n^{-1}y) = \|x+n^{-1}y\| \geq \|x\|.
\]
In particular, we have
\[
f_n(x) \geq \|x\|-\frac{1}{n}f_n(y) \geq \|x\|-\frac{1}{n}\|y\|
\]
and
\[
f_n(y) \geq n(\|x\|-f_n(x)) \geq 0
\]
for each $n \in \mathbb{N}$. Since $B_{X^*}$ is weakly$^*$ compact by the Banach-Alaoglu theorem, there exists a subnet $(f_{n_a})_a$ of $(f_n)_n$ that converges weakly$^*$ to some $f_0 \in B_{X^*}$. It follows that
\[
f_0(x)=\lim_a f_{n_a}(x)=\lim_n f_n (x) =\|x\| ,
\]
and $f_0(y) =\lim_a f_{n_a}(y) \geq 0$; that is, $f_0 \in \nu (x)$ and $f_0(y) \geq 0$.

Conversely, if $f(y) \geq 0$ for some $f \in \nu (x)$, then
\[
\|x+\lambda y\| \geq f(x+\lambda y) = \|x\|+\lambda f(y) \geq \|x\|
\]
for each $\lambda \geq 0$. Hence, $x\perp_{BJ}^+y$ holds.

(ii) If $x \perp_{BJ}^{++}y$, then $x\perp_{BJ}^+y$ and $x \not \perp_{BJ} y$. By (i), we obtain $f_1(y) \geq 0$ for some $f_1 \in \nu (x)$. Moreover, by Lemma~\ref{James}, $f(y) \neq 0$ for each $f \in \nu (x)$. It follows that $f_1(y)>0$. Now, suppose that $f_2(y)<0$ for some $f_2 \in \nu (x)$. Since $\nu (x)$ is convex, we derive
\[
g=\frac{-f_2(y)}{f_1(y)-f_2(y)}f_1+\frac{f_1(y)}{f_1(y)-f_2(y)}f_2 \in \nu (x),
\]
and $g(y)=0$. However, this contradicts $x \not \perp_{BJ}y$ by Lemma~\ref{James}. Thus, there is no $f \in \nu (x)$ satisfying $f(y)\leq 0$; that is, we have $f(y)>0$ for each $f \in \nu (x)$.

The converse follows from (i) and Lemma~\ref{James}.
\end{proof}
%%%%%%%%%%%%%%%%%%%%%%%%%%%%%%%%%%%%%%%%%%%%%%%%%%%%%%%%%%
In the following, let $R_x=\{ y \in X: x\perp_{BJ}y\}$ and $R_x^{++} = \{ y \in X: x \perp_{BJ}^{++} y\}$ for each $x \in X$. Below, we show that $R_x^{++}$ is an open convex subset of $X$ which is maximal as a connected subset of $X \setminus R_X$.
%%%%%%%%%%%%%%%%%%%%%%%%%%%%%%%%%%%%%%%%%%%%%%%%%%%%%%%%%%
\begin{lemma}\label{convex}

Let $X$ be a Banach space, and let $x \in X \setminus \{0\}$. Then, $R_x^{++}$ is an open convex subset of $X \setminus R_x$ containing $x$.
\begin{proof}
Set $\Pi^{++}(f)=\{ y \in X : f(y)>0\}$ for each $f \in X^*$. We then have $R_x^{++}=\bigcap \{\Pi^{++}(f) : f\in \nu (x)\}$ by Lemma~\ref{acute} (ii). Since $\Pi^{++}(f)$ is open, convex, and contains $x$ whenever $f \in \nu (x)$, it follows that $R_x^{++}$ is an open convex subset of $X \setminus R_x$ containing $x$.
\end{proof}
\end{lemma}
%%%%%%%%%%%%%%%%%%%%%%%%%%%%%%%%%%%%%%%%%%%%%%%%%%%%%%%%%%
\begin{lemma}\label{connected}

Let $X$ be a Banach space, and let $x \in X \setminus \{0\}$. Suppose that $A$ is a connected subset of $X \setminus R_x$ containing $x$. Then, $A \subset R_x^{++}$.
\begin{proof}
To show $A \subset R_x^{++}$, suppose to the contrary that $A \not \subset R_x^{++}$. Let $y \in A \setminus R_x^{++}$. Then, $f_0(y) \leq 0$ for some $f_0 \in \nu (x)$ by Lemma~\ref{acute}. However, since $A$ is connected and $f_0(x)=\|x\|>0$, the intermediate value theorem ensures that $f_0(z)=0$ for some $z \in A$. It follows from Lemma~\ref{James} that $x \perp_{BJ}z$. This contradicts $A \subset X \setminus R_x$. Therefore, $A \subset R_x^{++}$ holds.
\end{proof}
\end{lemma}
%%%%%%%%%%%%%%%%%%%%%%%%%%%%%%%%%%%%%%%%%%%%%%%%%%%%%%%%%%
Based on these lemmas, it turns out that if a Birkhoff-James orthogonality preserver is continuous, it then also preserves the Birkhoff-James acute angles.
%%%%%%%%%%%%%%%%%%%%%%%%%%%%%%%%%%%%%%%%%%%%%%%%%%%%%%%%%%
\begin{theorem}\label{preserver}

Let $X,Y$ be Banach spaces, and let $T:X \to Y$ be a Birkhoff-James orthogonality preserver. Suppose that $T$ is continuous. Then, $x,y \in X$ and $x \perp_{BJ}^+ y$ implies $Tx \perp_{BJ}^+ Ty$. Consequently, if $T$ is bicontinuous, then $T$ is a Birkhoff-James acute angle preserver.
\begin{proof}
We first note that $T(R_x)=R_{Tx}$ holds for each $x \in X$. Take an arbitrary $x \in X \setminus \{0\}$. Then, by Lemma~\ref{convex}, $R_x^{++}$ is an open convex subset of $X \setminus R_x$ containing $x$. Since $T$ is continuous, it follows that $T(R_x^{++})$ is a connected subset of $Y \setminus R_{Tx}$ containing $Tx$. Hence, by Lemma~\ref{connected}, we obtain $T(R_x^{++})\subset R_{Tx}^{++}$.

Now, suppose that $x,y \in X$, and $x \perp_{BJ}^+ y$. If $x \perp_{BJ} y$, then $Tx \perp_{BJ} Ty$, whereas if $x \not \perp_{BJ} y$, it then follows that $y \in R_x^{++}$, and $Ty \in R_{Tx}^{++}$. Therefore, in either case, we obtain $Tx \perp_{BJ}^+ Ty$.
\end{proof}
\end{theorem}
%%%%%%%%%%%%%%%%%%%%%%%%%%%%%%%%%%%%%%%%%%%%%%%%%%%%%%%%%%
Let $X,Y$ be Banach spaces, and let $T:X \to Y$. Then, $T$ is said to be \emph{norm-preserving} if $\|Tx\|=\|x\|$ for each $x \in X$, and \emph{homogeneous} if $T(cx)=cTx$ for each $x \in X$ and $c \in \mathbb{R}$. We now introduce a new notion of the equivalence of Banach spaces that is stronger than ``$\sim_{BJ}$.''
%%%%%%%%%%%%%%%%%%%%%%%%%%%%%%%%%%%%%%%%%%%%%%%%%%%%%%%%%%
\begin{definition}

Let $X,Y$ be Banach spaces. Then, $X$ is said to be \emph{strongly isomorphic} to $Y$ with respect to the structure of Birkhoff-James orthogonality, denoted by $X=_{BJ}Y$, if there exists a norm-preserving homogeneous bicontinuous Birkhoff-James orthogonality preserver $T:X \to Y$.
\end{definition}
%%%%%%%%%%%%%%%%%%%%%%%%%%%%%%%%%%%%%%%%%%%%%%%%%%%%%%%%%%
\begin{remark}

Obviously, $X=Y$ implies $X=_{BJ}Y$, and $X=_{BJ}Y$ implies $X\sim_{BJ}Y$.
\end{remark}
%%%%%%%%%%%%%%%%%%%%%%%%%%%%%%%%%%%%%%%%%%%%%%%%%%%%%%%%%%
The rest of this section is devoted to showing that if $X=_{BJ}Z$ and $Y=_{BJ}W$, then $X\oplus_\infty Y=_{BJ}Z\oplus_\infty W$, where $X \oplus_\infty Y$ is the Banach space $X\times Y$ endowed with the norm $\|(x,y)\|_\infty = \max \{\|x\|,\|y\|\}$. We begin with a further analysis of Birkhoff-James acute angles.
%%%%%%%%%%%%%%%%%%%%%%%%%%%%%%%%%%%%%%%%%%%%%%%%%%%%%%%%%%
\begin{lemma}\label{limit}

Let $X$ be a Banach space, and let $x,y \in X$. Then, $x \perp_{BJ}^+ y$ if and only if there exists a sequence of positive real numbers $(\lambda_n)_n$ such that $\lim_n \lambda_n =0$ and $\|x+\lambda_n y\| \geq \|x\|$ for each $n$.
\begin{proof}
It is sufficient to show the ``if'' part. Suppose that $(\lambda_n)_n$ is a sequence of positive real numbers such that $\lim_n \lambda_n =0$ and $\|x+\lambda_n y\| \geq \|x\|$ for each $n$. Then, for each $n$, there exists an $f_n \in S_{X^*}$ such that
\[
f_n(x+\lambda_n y)=\|x+\lambda_n y\| \geq \|x\|
\]
for each $n$. By the Banach-Alaoglu theorem, we have a subnet $(f_{n_a})_a$ of $(f_n)_n$ that converges weakly$^*$ to some $f_0 \in B_{X^*}$. Since
\[
f_n (x) \geq \|x\|-\lambda_n f_n(y) \geq \|x\|-\lambda_n \|y\|
\]
and
\[
f_n(y) \geq \lambda_n^{-1}(\|x\|-f_n(x)) \geq 0
\]
for each $n$, it follows that $f_0(x)=\|x\|$ and $f_0(y) \geq 0$. Therefore, $x \perp_{BJ}^+y$ by Lemma~\ref{acute} (i).
\end{proof}
\end{lemma}
%%%%%%%%%%%%%%%%%%%%%%%%%%%%%%%%%%%%%%%%%%%%%%%%%%%%%%%%%%
\begin{lemma}\label{either}

Let $X,Y$ be Banach spaces, and let $(x_1,y_1),(x_2,y_2) \in X \oplus_\infty Y$. Then, $(x_1,y_1) \perp_{BJ}^+ (x_2,y_2)$ if and only if one of the following holds:
\begin{itemize}
\item[{\rm (i)}] $\|x_1\|>\|y_1\|$, and $x_1 \perp_{BJ}^+ x_2$;
\item[{\rm (ii)}] $\|x_1\|=\|y_1\|$, and $x_1 \perp_{BJ}^+ x_2$ or $y_1 \perp_{BJ}^+ y_2$;
\item[{\rm (iii)}] $\|x_1\|<\|y_1\|$, and $y_1 \perp_{BJ}^+ y_2$.
\end{itemize}
\begin{proof}
Suppose that $(x_1,y_1) \perp_{BJ}^+ (x_2,y_2)$. Then, either $\|x_1\|>\|y_1\|$, $\|x_1\|=\|y_1\|$, or $\|x_1\|<\|y_1\|$ holds. If $\|x_1\|>\|y_1\|$, we have $\|(x_1,y_1)\|_\infty = \|x_1\|$. Moreover, $\|x_1+n^{-1}x_1\|>\|y_1+n^{-1}y_2\|$ for a sufficiently large $n$, and thus
\[
\|x_1+n^{-1}x_2\|=\|(x_1,y_1)+n^{-1}(x_2,y_2)\|_\infty \geq \|(x_1,y_1)\|_\infty=\|x_1\| .
\]
Hence, by Lemma~\ref{limit}, we obtain $x_1 \perp_{BJ}^+ x_2$. In the case of $\|x_1\|<\|y_1\|$, we know $y_1 \perp_{BJ}^+ y_2$ in a similar way.

Finally, we consider the case of $\|x_1\|=\|y_1\|$. Then, $\|(x_1,y_1)\|_\infty =\|x_1\|=\|y_1\|$. Set
\begin{align*}
N_1&=\{ n \in \mathbb{N} : \|x_1+n^{-1}x_2\| \geq \|y_1+n^{-1}y_2\|\}\\
N_2&=\{ n \in \mathbb{N} : \|x_1+n^{-1}x_2\| < \|y_1+n^{-1}y_2\|\} .
\end{align*}
Since $N_1 \cup N_2=\mathbb{N}$, either $N_1$ or $N_2$ is an infinite set. If $N_1$ is an infinite set, as in the preceding paragraph, we have $x_1 \perp_{BJ}^+ x_2$. Similarly, if $N_2$ is an infinite set, $y_1 \perp_{BJ}^+ y_2$ holds.

Conversely, suppose that (i), (ii), or (iii) holds. If (i) holds, then
\[
\|(x_1,y_1)+\lambda (x_2,y_2)\|_\infty \geq \|x_1+\lambda x_2\| \geq \|x_1\|=\|(x_1,y_1)\|_\infty
\]
for each $\lambda \geq 0$. Hence, $(x_1,y_1) \perp_{BJ}^+ (x_2,y_2)$. A similar argument shows the same conclusion in the case of (iii). Assume that (ii) holds. Then, either
\[
\|x_1+\lambda x_2\| \geq \|x_1\|=\|(x_1,y_1)\|_\infty
\]
for each $\lambda \geq 0$, or 
\[
\|y_1+\lambda y_2\| \geq \|y_1\|=\|(x_1,y_1)\|_\infty
\]
for each $\lambda \geq 0$. Thus, in either case, it follows that
\[
\|(x_1,y_1)+\lambda (x_2,y_2)\|_\infty \geq \|y_1\|=\|(x_1,y_1)\|_\infty
\]
for each $\lambda \geq 0$. This proves that $(x_1,y_1) \perp_{BJ}^+ (x_2,y_2)$.
\end{proof}
\end{lemma}
%%%%%%%%%%%%%%%%%%%%%%%%%%%%%%%%%%%%%%%%%%%%%%%%%%%%%%%%%%
Let $X,Y$ be Banach spaces, and let $T:X \to Y$. Then, $T$ is called a \emph{Birkhoff-James acute angle preserver} if it is bijective and $x \perp_{BJ}^+y$ if and only if $Tx\perp_{BJ}^+Ty$. For homogeneous Birkhoff-James acute angle preservers, we have the following result:
%%%%%%%%%%%%%%%%%%%%%%%%%%%%%%%%%%%%%%%%%%%%%%%%%%%%%%%%%%
\begin{proposition}\label{acute-ortho}

Let $X,Y$ be Banach spaces, and let $T:X \to Y$ be a homogeneous Birkhoff-James acute angle preserver. Then, $T$ is a Birkhoff-James orthogonality preserver.
\begin{proof}
Suppose that $x,y \in X$ and $x \perp_{BJ} y$. Then, $x \perp_{BJ}^+ y$ and $x \perp_{BJ}^+ (-y)$; in addition, by $T(-y)=-Ty$, this implies that $Tx \perp_{BJ}^+Ty$ and $Tx \perp_{BJ}^+ (-Ty)$. Hence, it follows that $Tx \perp_{BJ} Ty$. This completes the proof since the exact same argument is valid for $T^{-1}$. 
\end{proof}
\end{proposition}
%%%%%%%%%%%%%%%%%%%%%%%%%%%%%%%%%%%%%%%%%%%%%%%%%%%%%%%%%%
We are now ready to prove the following theorem.

%%%%%%%%%%%%%%%%%%%%%%%%%%%%%%%%%%%%%%%%%%%%%%%%%%%%%%%%%%
\begin{theorem}\label{induction}

Let $X,Y,Z,W$ be Banach spaces such that $X=_{BJ}Z$ and $Y=_{BJ}W$. Then, $X \oplus_\infty Y =_{BJ} Z \oplus_\infty W$.
\begin{proof}
Let $S:X \to Z$ and $T:Y \to W$ be norm-preserving homogeneous bicontinuous Birkhoff-James preservers, and let
\[
R(x,y) = (Sx,Ty)
\]
for each $(x,y) \in X\oplus_\infty \mathbb{R}$. Then, $R: X \oplus_\infty Y \to Z\oplus_\infty W$ is a norm-preserving homogeneous bicontinuous bijection because $R^{-1}$ is given by $R^{-1}(z,w)=(S^{-1}z,T^{-1}w)$ for each $(z,w) \in Z \oplus_\infty W$. Moreover, $R,T$ preserve the Birkhoff-James acute angle by Theorem~\ref{preserver}. Now, suppose that $(x_1,y_1),(x_2,y_2) \in X \oplus_\infty Y$, and $(x_1,y_1) \perp_{BJ}^+ (x_2,y_2)$. From Lemma~\ref{either}, one of the following holds:
\begin{itemize}
\item[{\rm (i)}] $\|x_1\|>\|y_1\|$, and $x_1 \perp_{BJ}^+ x_2$;
\item[{\rm (ii)}] $\|x_1\|=\|y_1\|$, and $x_1 \perp_{BJ}^+ x_2$ or $y_1 \perp_{BJ}^+ y_2$;
\item[{\rm (iii)}] $\|x_1\|<\|y_1\|$, and $y_1 \perp_{BJ}^+ y_2$.
\end{itemize}
If (i) holds, then $\|Sx_1\|=\|x_1\|>\|y_1\|=\|Ty_1\|$ and $Sx_1 \perp_{BJ}^+Sx_2$; hence, from Lemma~\ref{either} again, $R(x_1,y_1) \perp_{BJ}^+ R(x_2,y_2)$. Similarly, we obtain the same conclusion for the other cases. Thus, $R$ preserves the Birkhoff-James acute angle. Since the same is true for $R^{-1}$, by Proposition~\ref{acute-ortho}, we can conclude that $R$ is a Birkhoff-James orthogonality preserver. Therefore, $X \oplus_\infty Y =_{BJ} Z \oplus_\infty W$.
\end{proof}
\end{theorem}
%%%%%%%%%%%%%%%%%%%%%%%%%%%%%%%%%%%%%%%%%%%%%%%%%%%%%%%%%%
\begin{corollary}\label{Z-same}

Let $X,Y,Z$ be Banach spaces such that $X=_{BJ}Y$. Then, $X \oplus_\infty Z=_{BJ} Y\oplus_\infty Z$.
\end{corollary}
%%%%%%%%%%%%%%%%%%%%%%%%%%%%%%%%%%%%%%%%%%%%%%%%%%%%%%%%%%

%%%%%%%%%%%%%%%%%%%%%%%%%%%%%%%%%%%%%%%%%%%%%%%%%%%%%%%%%%
%%%%%%%%%%%%%%%%%%%%%%%%%%%%%%%%%%%%%%%%%%%%%%%%%%%%%%%%%%

%%%%%%%%%%%%%%%%%%%%%%%%%%%%%%%%%%%%%%%%%%%%%%%%%%%%%%%%%%
%%%%%%%%%%%%%%%%%%%%%%%%%%%%%%%%%%%%%%%%%%%%%%%%%%%%%%%%%%
\section{Examples}

In this section, we construct example pairs of real Banach spaces $(X,Y)$ of arbitrary dimensions such that $X=_{BJ}Y$ and $X \neq Y$. First, we improve \cite[Theorem 4.7]{Tan22} by slightly modifying its proof.
%%%%%%%%%%%%%%%%%%%%%%%%%%%%%%%%%%%%%%%%%%%%%%%%%%%%%%%%%%
\begin{theorem}\label{smooth-Radon}

Let $X$ be a smooth Radon plane. Then, $X=_{BJ}\ell_2^2$.
\end{theorem}
\begin{proof}
Let $x(\theta )=(\cos \theta ,\sin \theta )$ and $y(\theta)=\|x(\theta)\|_X^{-1}x(\theta)$ for each $\theta \in [0,2\pi]$. Then, as in the proof of \cite[Theorem 4.7]{Tan22}, there exists a bijection $\eta :[0,\pi/2] \to [\pi /2,\pi]$ such that $\eta (0)=\pi /2$, $\eta (\pi/2) = \pi$, and $y(\theta ) \perp_{BJ} y(\eta (\theta ))$ for each $\theta \in [0,\pi /2]$. We first show that this $\eta$ is automatically continuous. Let $\theta_0 \in [0,\pi /2]$. Suppose that $(\theta_n)_n$ is a sequence in $[0,\pi/2]$ that converges to $\theta_0$. By the definition of $\eta$, we have $y(\theta_n ) \perp_{BJ}y(\eta (\theta_n))$, which implies that
\[
\|y(\theta_n) +\lambda y(\eta (\theta_n))\| \geq \|y(\theta_n)\|=1
\]
for each $\lambda \in \mathbb{R}$. Now, take an arbitrary subsequence $(\eta (\theta_{k_n}))_n$ of $(\eta (\theta_n))_n$. Then, there exists a subsequence $(\eta (\theta_{l_{k_n}}))_n$ of $(\eta (\theta_{k_n}))_n$ that converges to some $\theta_0' \in [\pi/2,\pi]$. Since the mapping $\theta \mapsto y(\theta )$ is a homeomorphism from $[0,\pi]$ to $S_X^+$, we derive
\[
\|y(\theta_0) +\lambda y(\theta_0')\|=\lim_n \|y(\theta_{l_{k_n}}) +\lambda y(\eta (\theta_{l_{k_n}}))\| \geq 1
\]
for each $\lambda \in \mathbb{R}$. This means that $y(\theta_0) \perp_{BJ}y(\theta_0')$; hence, $\theta_0'=\eta (\theta_0)$ by the smoothness of $X$. Therefore, $\eta (\theta_n) \to \eta (\theta_0)$; that is, $\eta$ is continuous on $[0,\pi/2]$.

Next, set
\[
T(x(\theta )) = \left\{ \begin{array}{ll}
y(\theta ) & (\theta \in [0,\pi/2])\\
y(\eta (\theta-\pi/2)) & (\theta \in [\pi /2,\pi]) 
\end{array}
\right. .
\]
Then, as was shown in \cite[Theorem 4.7]{Tan22}, $T:S_{\ell_2^2}^+\to S_X$ is a bijection such that $x(\theta ) \perp x(\theta')$ in $\ell_2^2$ if and only if $Tx(\theta ) \perp_{BJ}Tx(\theta' )$ in $X$. We note that the mapping $\theta \mapsto x(\theta )$ is a homeomorphism from $[0,\pi/2]$ onto $S_{\ell_2^2}^+$. Since $y(\pi/2)=y(\eta (0))$, and $\eta$ is continuous as indicated in the preceding paragraph, $T$ is also continuous.

We extend $T$ by setting $T^=(x)=-T(-x)$ for each $-S_{\ell_2^2}^+$. Then, $T^=:S_{\ell_2^2}\to S_X$ is a bijection such that $x \perp_{BJ}y$ in $\ell_2^2$ if and only if $T^=x\perp_{BJ}T^=y$. It is easy to check whether $T^=$ is continuous. In fact, $T^=$ is a homeomorphism since it is a continuous bijection from a compact space $S_{\ell_2^2}$ onto a Hausdorff space $S_X$.

Finally, we define the full extension $T^\sim$ of $T^=$ to $\ell_2^2$ by
\[
T^\sim x = \left\{ \begin{array}{ll}
0 & (x=0)\\
\displaystyle \|x\|_2T^=\left(\frac{1}{\|x\|_2}x\right) & (x \neq 0) 
\end{array}
\right. .
\]
This gives rise to a norm-preserving homogeneous continuous Birkhoff-James orthogonality preserver from $\ell_2^2$ onto $X$. Moreover, since
\[
(T^\sim)^{-1}y = \left\{ \begin{array}{ll}
0 & (y=0)\\
\displaystyle \|y\|_X(T^=)^{-1}\left(\frac{1}{\|y\|_X}y\right) & (y \neq 0) 
\end{array}
\right. ,
\]
$(T^\sim)^{-1}$ is also continuous. Thus, $X=_{BJ}\ell_2^2$ holds.
\end{proof}
%%%%%%%%%%%%%%%%%%%%%%%%%%%%%%%%%%%%%%%%%%%%%%%%%%%%%%%%%%
For a nonempty set $I$, the symbol $c_0(I)$ denotes the Banach space of all systems $(a_n)_{n \in I}$ such that $\{ n \in I: |a_n| \geq \varepsilon \}$ is finite for each $\varepsilon >0$. To prove the non-isometric clause of the examples, we need the following technical lemma.
%%%%%%%%%%%%%%%%%%%%%%%%%%%%%%%%%%%%%%%%%%%%%%%%%%%%%%%%%%
\begin{lemma}\label{non-Hilbert}

Let $X$ be a Banach space, and let $I$ be a nonempty set. If $X$ does not contain a two-dimensional Hilbert subspace, the same is then true for $X \oplus_\infty c_0(I)$.
\begin{proof}
We first prove the theorem under the case in which $I$ is a singleton. Suppose that $X \oplus_\infty \mathbb{R}$ contains a two-dimensional Hilbert subspace $M$. Let $\{(x,r),(y,s)\} \subset M$ be an orthonormal basis for $M$. We then have
\[
\|\alpha (x,r)+\beta (y,s)\|_\infty = (\alpha^2+\beta^2)^{1/2}
\] 
for each $\alpha ,\beta \in \mathbb{R}$. In particular, setting $(\alpha ,\beta)=(r,s)$, we obtain
\[
r^2+s^2\leq \|r (x,r)+s (y,s)\|_\infty = (r^2+s^2)^{1/2},
\]
which implies that $(r^2+s^2)^{1/2} \leq 1$. Combining this with the Cauchy-Schwartz inequality yields
\[
|\alpha r+\beta s| \leq (\alpha^2+\beta^2)^{1/2}
\]
for each $(\alpha ,\beta)$. Moreover, the equality holds only if $(\alpha ,\beta)=c(r,s)$ for some $c \in \mathbb{R}$. This means that
\[
\|\alpha x+\beta y\|=(\alpha^2+\beta^2)^{1/2}
\]
for each $(\alpha ,\beta) \in \mathbb{R}^2 \setminus [(r,s)]$. Since the set $\mathbb{R}^2\setminus [(r,s)]$ is dense in $\mathbb{R}^2$ with respect to the product topology on $\mathbb{R}^2$, it follows that
\[
\|\alpha x+\beta y\|=(\alpha^2+\beta^2)^{1/2}
\]
for all $(\alpha ,\beta) \in \mathbb{R}^2$. Therefore, $[x,y]$ is a two-dimensional Hilbert subspace of $X$.

Next, suppose that $I$ is finite, and that $X \oplus_\infty c_0(I)$ contains a two-dimensional Hilbert subspace. We may write $I=\{1,\ldots ,n\}$. Since
\[
X\oplus_\infty c_0(\{1,\ldots ,n\})=(X \oplus_\infty c_0(\{1,\ldots ,n-1\})) \oplus_\infty \mathbb{R} ,
\]
as indicated in the first paragraph of the proof, $X \oplus_\infty c_0(\{1,\ldots ,n-1\})$ contains a two-dimensional Hilbert subspace. Continuing this process, we infer that $X$ contains a two-dimensional Hilbert subspace.

Finally, assume that $I$ is infinite. Suppose that $\{(x,(a_n)_{n \in I}),(y,(b_n)_{n \in I})\} \subset X \oplus_\infty c_0(I)$ is an orthonormal basis for two-dimensional Hilbert subspace. Since $(a_n)_{n\in I},(b_n)_{n\in I} \in c_0(I)$, the sets
\begin{align*}
I_1=\left\{ n \in I : |a_n| \geq \frac{1}{\sqrt{2}} \right\} \\
I_2=\left\{ n \in I : |b_n| \geq \frac{1}{\sqrt{2}} \right\}
\end{align*}
are both finite. Set $I_0=I_1 \cup I_2$. Then, $I_0$ is also finite, and
\[
|\alpha a_n+\beta b_n| \leq (\alpha^2+\beta^2)^{1/2}(a_n^2+b_n^2)^{1/2}<(\alpha^2+\beta^2)^{1/2}
\]
for each $n \not \in I_0$ and each nonzero $(\alpha ,\beta) \in \mathbb{R}^2$. Hence, it follows that
\begin{align*}
(\alpha^2+\beta^2)^{1/2}
&=\|\alpha (x,(a_n)_{n\in I})+\beta (y,(b_n)_{n \in I})\|_\infty \\
&= \max \left\{\|\alpha x+\beta y\|,\max_{n \in I}|\alpha a_n+\beta b_n|\right\}\\
&= \max \left\{\|\alpha x+\beta y\|,\max_{n \in I_0}|\alpha a_n+\beta b_n|\right\}
\end{align*}
for each $(\alpha ,\beta ) \in \mathbb{R}^2$. This means that $X\oplus_\infty c_0(I_0)$ contains a two-dimensional Hilbert subspace. Thus, as described in the preceding paragraph, we can conclude that $X$ contains a two-dimensional Hilbert subspace.
\end{proof}
\end{lemma}
%%%%%%%%%%%%%%%%%%%%%%%%%%%%%%%%%%%%%%%%%%%%%%%%%%%%%%%%%%
Now, we prove the second main theorem of this paper.

%%%%%%%%%%%%%%%%%%%%%%%%%%%%%%%%%%%%%%%%%%%%%%%%%%%%%%%%%%
\begin{theorem}\label{main}

Let $X$ be a smooth Radon plane, and let $I$ be a nonempty set. Suppose that $X \neq \ell_2^2$. Then, $X \oplus_\infty c_0(I) =_{BJ}\ell_2^2 \oplus_\infty c_0(I)$, and $X \oplus_\infty c_0(I) \neq \ell_2^2 \oplus_\infty c_0(I)$.
\begin{proof}
By Corollary~\ref{Z-same} and Theorem~\ref{smooth-Radon}, we have $X \oplus_\infty c_0(I) =_{BJ}\ell_2^2 \oplus_\infty c_0(I)$. Moreover, since $X \neq \ell_2^2$, by Lemma~\ref{non-Hilbert}, $X \oplus_\infty c_0(I)$ does not contain two-dimensional Hilbert subspace, whereas $\ell_2^2 \subset \ell_2^2 \oplus_\infty c_0(I)$. Therefore, $X \oplus_\infty c_0(I) \neq \ell_2^2 \oplus_\infty c_0(I)$.
\end{proof}
\end{theorem}
%%%%%%%%%%%%%%%%%%%%%%%%%%%%%%%%%%%%%%%%%%%%%%%%%%%%%%%%%%
\begin{example}

As was shown in \cite[Example 4.8]{Tan22}, the Day-James space $\ell_{p,q}^2$ with $p^{-1}+q^{-1}=1$ is a non-Hilbert smooth Radon plane, where $\ell_{p,q}^2$ is the Banach space $\mathbb{R}^2$ endowed with the norm
\[
\|(a,b)\|_{p,q}=\left\{ \begin{array}{ll}
(|a|^p+|b|^p)^{1/p} & (ab \geq 0)\\
(|a|^q+|b|^q)^{1/q} & (ab \leq 0)
\end{array}
\right. .
\]
Hence, we have an infinitely many examples of $X$ that satisfy the assumption of the preceding theorem.
\end{example}
%%%%%%%%%%%%%%%%%%%%%%%%%%%%%%%%%%%%%%%%%%%%%%%%%%%%%%%%%%

%%%%%%%%%%%%%%%%%%%%%%%%%%%%%%%%%%%%%%%%%%%%%%%%%%%%%%%%%%
%%%%%%%%%%%%%%%%%%%%%%%%%%%%%%%%%%%%%%%%%%%%%%%%%%%%%%%%%%

%%%%%%%%%%%%%%%%%%%%%%%%%%%%%%%%%%%%%%%%%%%%%%%%%%%%%%%%%%
%%%%%%%%%%%%%%%%%%%%%%%%%%%%%%%%%%%%%%%%%%%%%%%%%%%%%%%%%%
\section{Remarks}

The proofs of two main theorems of this paper are only valid for real Banach spaces. Indeed, in a complex case, the operator $S$ constructed in the proof of Theorem~\ref{sm-eq} is generally \emph{semilinear}. To obtain the linearity or antilinearity of $S$, we need some arguments related to the continuity. For example, it was shown by \cite[Lemma 2 and Corollary]{KM46} that a semilinear bijection between infinite-dimensional complex Banach spaces is linear or antilinear if it preserves closed hyperplanes in both directions. Moreover, such an operator is automatically continuous; see Fillmore and Longstaff~\cite[Lemma 2 and Lemma 3]{FL86}. The finite-dimensional case is known to be more complicated. At this stage, we do not know if a Birkhoff-James orthogonality preserver between smooth complex Banach spaces satisfies some continuity conditions. At a minimum, if the considered complex Banach spaces are reflexive, smooth, and infinite-dimensional, we may thus overcome this difficulty by using \cite[Corollary 3.10]{Tan22} and \cite[Theorem 1]{FL86} instead of Theorem~\ref{sm-ref} and Lemma~\ref{Mackey}. However, there is no idea on how to fill in the gaps if the spaces are non-reflexive, or finite-dimensional.

The validity of the second main theorem for complex spaces is under a more serious condition. The construction of examples described herein relies heavily on the existence of non-Hilbert Radon planes. However, to the best of the author's knowledge, it is not known if a complex Radon plane exists.

Finally, we note that all known pairs of real Banach spaces $(X,Y)$ satisfying $X \sim_{BJ}Y$ are mutually isomorphic. Hence, we wonder if $X \sim_{BJ}Y$ implies $X \cong Y$ without any assumption. Meanwhile, by \cite[Remarks 3.14 and 3.15]{Tan22}, we have a complex Banach space $X$ that satisfies (automatically) $X\sim_{BJ}\overline{X}$, and $X \not \cong \overline{X}$. Therefore, the corresponding problem in the complex case would be the following: Does $X \sim_{BJ}Y$ imply $X \cong Y$ or $X\cong \overline{Y}$?

We conclude this paper with some useful tools for analyzing the structure of Birkhoff-James orthogonality in Banach spaces. In 2005, Turn\v{s}ek~\cite{Tur05} considered a local symmetry condition to Birkhoff-James orthogonality, with was then applied to determine the forms of the Birkhoff-James orthogonality preservers on $B(H)$, where $B(H)$ is the Banach space of all bounded linear operators on a complex Hilbert space $H$. It was shown accurately in~\cite[Theorem 2.5]{Tur05} that if $B \in B(H)$, and if $A \perp_{BJ}B$ always implies $B \perp_{BJ}A$ in $B(H)$, then $B$ is a scalar multiple of an isometry or a coisometry. Moreover, this result was generalized to a characterization of extreme points of $B_A$ for a von Neumann algebra $A$; see~\cite[Theorem 4.7]{KST18}. Now, such a property of $B$ is known as \emph{right symmetry} for Birkhoff-James orthogonality. Similarly, an element $x$ of a Banach space $X$ is called a \emph{left symmetric point} for Birkhoff-James orthogonality if $x \perp_{BJ}y$ always implies $y \perp_{BJ}x$ in $X$. The terms ``left symmetric'' and ``right symmetric'' were introduced by Sain~\cite{Sai17}, and Sain, Ghosh, and Paul~\cite{SGP17}, respectively. These locally symmetric points for Birkhoff-James orthogonality reflect the geometric features of Banach spaces, and are preserved under Birkhoff-James orthogonality preservers. Hence, they may play important roles in classifying general Banach spaces based on the structure of Birkhoff-James orthogonality.

Another interesting tool is \emph{orthographs} of \emph{symmetrized} Birkhoff-James orthogonality that were recently introduced and studied by Aramba\v{s}i\'{c} et al.~\cite{AGKRZ21}. We define $x\perp y$ by $x\perp_{BJ}y$ and $y \perp_{BJ}x$. Then, ``$\perp$'' is called \emph{mutual} Birkhoff-James orthogonality. The orthograph $\Gamma (X)$ of a Banach space $X$ induced by ``$\perp$'' is defined as follows: The vertex set $V(\Gamma (X))$ is the set of all one-dimensional subspaces of $X$, and the vertices $[x],[y]$ are \emph{adjacent} if $x\perp y$. Naturally, the properties of the orthographs depend on the geometric structure of the Banach spaces. The author believes that orthographs are preserved under ``$\sim_{BJ}$'' in a reasonable manner, and are useful in some important circumstances.
%%%%%%%%%%%%%%%%%%%%%%%%%%%%%%%%%%%%%%%%%%%%%%%%%%%%%%%%%%
%%%%%%%%%%%%%%%%%%%%%%%%%%%%%%%%%%%%%%%%%%%%%%%%%%%%%%%%%%
\section*{Acknowledgment}

The author would like to thank Editage (www.editage.com) for English language editing. This work was supported by JSPS KAKENHI Grant Number JP19K14561.

%%%%%%%%%%%%%%%%%%%%%%%%%%%%%%%%%%%%%%%%%%%%%%%%%%%%%%%%%%
%%%%%%%%%%%%%%%%%%%%%%%%%%%%%%%%%%%%%%%%%%%%%%%%%%%%%%%%%%

%%%%%%%%%%%%%%%%%%%%%%%%%%%%%%%%%%%%%%%%%%%%%%%%%%%%%%%%%%
\end{document}